# When is Eaton's Markov chain irreducible?

JAMES P. HOBERT[*], AIXIN TAN and RUITAO LIU

*Department of Statistics, University of Florida, Gainesville, FL 32611, USA.*
*E-mail:* [*]*jhobert@stat.ufl.edu*

Consider a parametric statistical model $P(\mathrm{d}x|\theta)$ and an improper prior distribution $\nu(\mathrm{d}\theta)$ that together yield a (proper) formal posterior distribution $Q(\mathrm{d}\theta|x)$. The prior is called *strongly admissible* if the generalized Bayes estimator of every bounded function of $\theta$ is admissible under squared error loss. Eaton [*Ann. Statist.* **20** (1992) 1147–1179] has shown that a sufficient condition for strong admissibility of $\nu$ is the *local recurrence* of the Markov chain whose transition function is $R(\theta, \mathrm{d}\eta) = \int Q(\mathrm{d}\eta|x) P(\mathrm{d}x|\theta)$. Applications of this result and its extensions are often greatly simplified when the Markov chain associated with $R$ is irreducible. However, establishing irreducibility can be difficult. In this paper, we provide a characterization of irreducibility for general state space Markov chains and use this characterization to develop an easily checked, necessary and sufficient condition for irreducibility of Eaton's Markov chain. All that is required to check this condition is a simple examination of $P$ and $\nu$. Application of the main result is illustrated using two examples.

## 1. Introduction

Consider a parametric statistical decision problem with sample space $\mathsf{X}$ and parameter space $\Theta$. Both $\mathsf{X}$ and $\Theta$ are assumed to be Polish spaces equipped with their Borel $\sigma$-algebras $\mathcal{B}(\mathsf{X})$ and $\mathcal{B}(\Theta)$. Suppose that $P : \mathcal{B}(\mathsf{X}) \times \Theta \to [0,1]$ represents a parametric statistical model, that is, for each $\theta$, $P(\cdot|\theta)$ is a probability measure and, for each $A$, $P(A|\cdot)$ is a measurable function. As usual, the idea is that we will observe a random element whose distribution is $P(\mathrm{d}x|\theta)$, the goal being to use the observation to make inferences about the unknown parameter $\theta$. This will be done within the Bayesian paradigm using an *improper* prior distribution. In particular, let $\nu(\mathrm{d}\theta)$ denote a $\sigma$-finite measure with $\nu(\Theta) = \infty$. Define the marginal measure as

$$M(\mathrm{d}x) = \int_\Theta P(\mathrm{d}x|\theta) \nu(\mathrm{d}\theta).$$

Eaton [2] shows that if $M$ is $\sigma$-finite, then there exists a formal posterior distribution $Q$, defined as follows.

**Definition 1.** *A function $Q : \mathcal{B}(\Theta) \times \mathsf{X} \to [0,1]$ is called a* formal posterior distribution (FPD) *if:*

1. *$Q(\cdot|x)$ is a probability measure for each $x$;*
2. *$Q(B|\cdot)$ is a measurable function for each $B$;*
3. *$Q(\mathrm{d}\theta|x)M(\mathrm{d}x) = P(\mathrm{d}x|\theta)\nu(\mathrm{d}\theta)$, that is, for all $A \in \mathcal{B}(\mathsf{X})$ and $B \in \mathcal{B}(\Theta)$,*

$$\int_A Q(B|x)M(\mathrm{d}x) = \int_B P(A|\theta)\nu(\mathrm{d}\theta).$$

The FPD is unique in the sense that if $\tilde{Q}$ is another FPD, then there is an $M$-null set $A_0$ such that $x \notin A_0$ implies $Q(\mathrm{d}\theta|x) = \tilde{Q}(\mathrm{d}\theta|x)$. Throughout this paper, $M$ is assumed to be $\sigma$-finite, so an FPD is guaranteed to exist. We now briefly describe a method of evaluating improper prior distributions that is due to M.L. Eaton. (For a more in depth review of this area, see Eaton [4, 6].)

Consider the problem of estimating a bounded, real-valued function $\gamma(\theta)$ under squared error loss. Of course, the formal Bayes estimator of $\gamma(\theta)$ is $\hat{\gamma}(x) = \int_\Theta \gamma(\theta) Q(\mathrm{d}\theta|x)$. The risk function of a generic estimator, say $\delta$, is its mean squared error, that is,

$$r(\delta, \theta) = \int_\mathsf{X} (\delta(x) - \gamma(\theta))^2 P(\mathrm{d}x|\theta).$$

The estimator $\delta$ is called *almost-$\nu$-admissible* if for any estimator $\delta'$ such that

$$r(\delta', \theta) \leq r(\delta, \theta) \qquad \forall \theta \in \Theta,$$

the set $\{\theta \in \Theta : r(\delta', \theta) < r(\delta, \theta)\}$ has $\nu$-measure zero. If $P(\mathrm{d}x|\theta)$ and $\nu(\mathrm{d}\theta)$ combine to yield an FPD that generates (almost) admissible estimators for a large class of functions of $\theta$, then we would be willing to endorse $\nu$ as a good "all purpose" prior to use in conjunction with the statistical model $P(\mathrm{d}x|\theta)$. This idea provides motivation for the following definition.

**Definition 2.** *The prior $\nu$ is called* strongly admissible *if $\hat{\gamma}$ is almost-$\nu$-admissible for every bounded, real-valued function $\gamma$.*

Eaton [3] developed a sufficient condition for strong admissibility that involves the Markov transition function $R : \Theta \times \mathcal{B}(\Theta) \to [0,1]$ given by

$$R(\theta, \mathrm{d}\eta) = \int_\mathsf{X} Q(\mathrm{d}\eta|x) P(\mathrm{d}x|\theta).$$

Before we can state the result, we need a couple of concepts from general state space Markov chain theory. Let $W = \{W_n\}_{n=0}^\infty$ denote the Markov chain on $\Theta$ driven by $R$ and let $\Pr_\theta$ denote the overall probability law governing the chain when $W_0 = \theta$. For $B \in \mathcal{B}(\Theta)$, let $\sigma_B$ denote the first return to $B$, that is,

$$\sigma_B = \min\{n \geq 1 : W_n \in B\},$$

with the understanding that $\sigma_B = \infty$ if $W_n \in \overline{B}$ for all $n \geq 1$, where $\overline{B}$ denotes the complement of $B$.

**Definition 3.** *The Markov chain $W$ is called* locally-$\nu$-recurrent *if, for each $B$ with $0 < \nu(B) < \infty$, the set*

$$\{\theta \in B : \Pr{}_\theta(\sigma_B < \infty) < 1\}$$

*has $\nu$-measure 0.*

In words, the chain is locally-$\nu$-recurrent if, when started inside the set $B$, aside from a set of starting values that has $\nu$-measure 0, the chain returns to $B$ with probability 1. Note that, unlike the standard definition of recurrence (see, e.g., Meyn and Tweedie [11], Chapter 8), this definition pertains to both reducible and irreducible chains. Indeed, just before defining local-$\nu$-recurrence on page 1174, Eaton [3] states: "The following definition, a modified notion of recurrence, allows us to circumvent a discussion of irreducibility issues while relating our previous admissibility results to the recurrence of $W$." The following was proven in Eaton [3] (see also Eaton [5]).

**Theorem 1.** *If $W$ is a locally-$\nu$-recurrent Markov chain, then $\nu$ is a strongly admissible prior.*

Establishing local-$\nu$-recurrence directly using the definition (or the characterization based on the Dirichlet form of $R$) is typically infeasible. However, if $W$ is $\nu$-irreducible, that is, any set $B$ with $\nu(B) > 0$ is accessible from any $\theta \in \Theta$, then recurrence and local-$\nu$-recurrence are equivalent (Eaton, Hobert and Jones [7]). Hence, if $\nu$-irreducibility of $W$ can be demonstrated, then all of the techniques that have been developed for establishing recurrence can be brought to bear on the problem. Indeed, nearly all of the applications of Theorem 1 have involved first demonstrating that $W$ is $\nu$-irreducible and then showing that $W$ is recurrent. Examples can be found in Eaton [3], Hobert and Robert [9], Hobert and Schweinsberg [10] and Hobert, Marchev and Schweinsberg [8]. Similarly, Eaton *et al.* [7] have recently extended and generalized the theoretical results of Eaton [3] and Hobert and Robert [9] under the assumption that the chains of interest are irreducible.

There is one very simple sufficient condition for $\nu$-irreducibility of $W$ and this was used in most of the applications mentioned above. If the support of the statistical model does not depend on the parameter, that is, if the set $\{A \in \mathcal{B}(\mathsf{X}) : P(A|\theta) > 0\}$ is the same for all $\theta \in \Theta$, then $W$ is $\nu$-irreducible (Eaton *et al.* [7]). Until now, however, there has been no easy way to check for $\nu$-irreducibility of Eaton's Markov chain when this condition fails. In this paper, we provide an easily checked, necessary and sufficient condition for $\nu$-irreducibility of $W$. This result cannot be stated precisely at this point, but the sufficiency half, which is the practically important part, can be.

**Theorem 2.** *The Markov chain $W$ is $\nu$-irreducible if there do not exist two sets $A \in \mathcal{B}(\mathsf{X})$ and $C \in \mathcal{B}(\Theta)$ with the following properties: $C$ is non-empty, $\nu(\overline{C}) > 0$, $P(\overline{A}|\theta) = 0$ for every $\theta \in C$ and $P(A|\theta) = 0$ for $\nu$-almost all $\theta \in \overline{C}$.*

This result allows one to establish $\nu$-irreducibility of Eaton's Markov chain through a simple examination of $P$ and $\nu$. Neither the posterior distribution $Q$ nor the Markov transition function $R$ is required to check the condition. It is interesting to note that if the sets $A$ and $C$ do exist, then $P(\mathrm{d}x|\theta)$ and $P(\mathrm{d}x|\theta')$ are mutually singular probability measures whenever $\theta \in C$ and $\theta' \in \overline{C}$ (aside from a $\nu$-null set of $\theta$'s in $\overline{C}$). Thus, the statistical model is, in a sense, an artificial concatenation of two different models.

The remainder of the paper is organized as follows. Section 2 contains a new characterization of irreducibility for general state space Markov chains. This characterization is used in Section 3 to prove the main result. Application of the main result is illustrated in Section 4.

## 2. A characterization of irreducibility for general Markov chains

Let $S(y, \mathrm{d}z)$ be a Markov transition function on a general state space $(\mathsf{Y}, \mathcal{B}(\mathsf{Y}))$, as described, for example, in Meyn and Tweedie [11], Section 3.4. Denote the corresponding Markov chain by $Y = \{Y_n\}_{n=0}^{\infty}$. For $n \in \mathbb{N} := \{1, 2, 3, \ldots\}$, let $S^n(y, \mathrm{d}z)$ denote the $n$-step Markov transition function corresponding to $S$, which is defined inductively by

$$S^{n+1}(y, \mathrm{d}z) = \int_{\mathsf{Y}} S^n(w, \mathrm{d}z) S(y, \mathrm{d}w),$$

where $S^1 \equiv S$. Of course, $S^n(y, A) = \mathrm{Pr}_y(Y_n \in A)$, where $\mathrm{Pr}_y(\cdot)$ denotes the overall law governing $Y$ on $\mathsf{Y}^{\infty}$, assuming that $Y_0 = y$. Let $\varphi$ denote a non-trivial, $\sigma$-finite measure on $(\mathsf{Y}, \mathcal{B}(\mathsf{Y}))$. The following is a standard definition of irreducibility for general state space Markov chains.

**Definition 4** (Meyn and Tweedie [11], page 87). *The Markov chain $Y$ is called $\varphi$-irreducible if, for every measurable $A$ with $\varphi(A) > 0$ and every $y \in \mathsf{Y}$, there exists an $n \in \mathbb{N}$ (which may depend on $y$ and $A$) such that $S^n(y, A) > 0$.*

In words, $\varphi$-irreducibility means that every set $A$ with $\varphi(A) > 0$ is *accessible* from any $y \in \mathsf{Y}$. We will call the Markov chain $\varphi$-*reducible* when it is not $\varphi$-irreducible, that is, when there exist $y$ and $A$ with $\varphi(A) > 0$ such that $S^n(y, A) = 0$ for all $n \in \mathbb{N}$. We will sometimes find it convenient to apply the descriptions "$\varphi$-irreducible" and "$\varphi$-reducible" to the Markov transition function $S$. Our first result shows that if the chain is $\varphi$-reducible, then we may assume that $y \in \overline{A}$.

**Proposition 1.** *The Markov chain $Y$ is $\varphi$-reducible if and only if there exists $A \in \mathcal{B}(\mathsf{Y})$ with $\varphi(A) > 0$ and $y \in \overline{A}$ such that $S^n(y, A) = 0$ for all $n \in \mathbb{N}$.*

**Proof.** The sufficiency part is obvious. Now, assume that the chain is $\varphi$-reducible, so there exist $y \in \mathsf{Y}$ and $A \in \mathcal{B}(\mathsf{Y})$ with $\varphi(A) > 0$ such that $S^n(y, A) = 0$ for all $n \in \mathbb{N}$.

If $y \in \overline{A}$, then there is nothing to prove, so assume that $y \in A$. We will establish the existence of $y' \in \overline{A}$ such that $S^n(y', A) = 0$ for all $n \in \mathbb{N}$. First, for each $m \in \mathbb{N}$, define

$$B_m = \{w \in \overline{A} : S^m(w, A) > 0\},$$

and set $B = \bigcup_{m=1}^{\infty} B_m$. Now, fix $m \in \mathbb{N}$ and note that

$$0 = S^{m+1}(y, A) = \int_Y S^m(z, A) S(y, \mathrm{d}z) \geq \int_{B_m} S^m(z, A) S(y, \mathrm{d}z).$$

Since $S^m(z, A) > 0$ for $z \in B_m$, we must have $S(y, B_m) = 0$. But this result holds for every $m \in \mathbb{N}$, so it follows that $S(y, B) = 0$. Note that $\mathsf{Y}$ can be partitioned into $A$, $B$ and $\overline{A} \setminus B$ and we know that $S(y, A) = S(y, B) = 0$. Therefore, $S(y, \overline{A} \setminus B) = 1$, which implies that $\overline{A} \setminus B$ is not empty. Clearly, any $y' \in \overline{A} \setminus B$ satisfies $S^n(y', A) = 0$ for all $n \in \mathbb{N}$. □

In the classical case where $\mathsf{Y}$ is countable, the chain is called *irreducible* (with no prefix) if, for each $i, j \in \mathsf{Y}$, there exists an $n \in \mathbb{N}$ such that $S^n(i, \{j\}) > 0$. This is equivalent to $c$-irreducibility, where $c$ denotes counting measure on $\mathsf{Y}$. In this context, a non-empty set $C \subset \mathsf{Y}$ is called *closed* if, once the chain enters $C$, it cannot leave. Formally, $C$ is closed if $\sum_{j \in C} S(i, \{j\}) = 1$ for all $i \in C$. Obviously, the state space $\mathsf{Y}$ is closed. In fact, the Markov chain is irreducible if and only if $\mathsf{Y}$ has no proper, closed subset (see, e.g., Billingsley [1], Problem 8.21). We now extend these ideas to handle Markov chains on general state spaces.

**Definition 5.** *A set $C \in \mathcal{B}(\mathsf{Y})$ is called* closed *if it is non-empty and $S(y, \overline{C}) = 0$ for all $y \in C$.*

The following is the general state space version of the result in Billingsley's [1] Problem 8.21.

**Theorem 3.** *The Markov chain $Y$ is $\varphi$-reducible if and only if there exists a closed set $C$ with $\varphi(\overline{C}) > 0$.*

**Proof.** To prove sufficiency, suppose that $C$ is a closed set with $\varphi(\overline{C}) > 0$. Assume that for some $n \in \mathbb{N}$, $S^n(y, \overline{C}) = 0$ for all $y \in C$. Then, for any $y \in C$, we have

$$S^{n+1}(y, \overline{C}) = \int_Y S^n(z, \overline{C}) S(y, \mathrm{d}z) = \int_C S^n(z, \overline{C}) S(y, \mathrm{d}z) = 0.$$

Hence, by induction, $S^n(y, \overline{C}) = 0$ for all $y \in C$ and all $n \in \mathbb{N}$. Therefore, since $C$ is non-empty and $\varphi(\overline{C}) > 0$, the chain is $\varphi$-reducible.

Now, to prove necessity, assume that the chain is $\varphi$-reducible. By Proposition 1, there exist a measurable $A$ with $\varphi(A) > 0$ and a $y \in \overline{A}$ such that $S^n(y, A) = 0$ for all $n \in \mathbb{N}$. For each $m \in \mathbb{N}$, define

$$B_m = \{w \in \mathsf{Y} : S^m(w, A) > 0\}$$

and set $B = \bigcup_{m=1}^{\infty} B_m$. We will show that the measurable set $C := \overline{A} \cap \overline{B}$ is the closed set that we seek. First, $y \in C$, so $C$ is non-empty. Now, suppose that $y' \in C$ is such that $S(y', \overline{C}) > 0$. Then,

$$S(y', A) + S(y', B) \geq S(y', A \cup B) = S(y', \overline{C}) > 0.$$

Since $y' \in C \Rightarrow y' \in \overline{B} \Rightarrow y' \in \overline{B}_1$, we know that $S(y', A) = 0$, so it must be the case that $S(y', B) > 0$. Hence, there must exist a $k \in \mathbb{N}$ such that $S(y', B_k) > 0$, which implies that

$$S^{k+1}(y', A) = \int_{\mathsf{Y}} S^k(z, A) S(y', \mathrm{d}z) = \int_{B_k} S^k(z, A) S(y', \mathrm{d}z) > 0.$$

But this implies that $y' \in B_{k+1}$, which contradicts the fact that $y' \in C$. Therefore, it must be the case that $S(y, \overline{C}) = 0$ for all $y \in C$ and this implies that $C$ is closed. Finally, $\varphi(\overline{C}) \geq \varphi(A) > 0$. $\square$

In the next section, we use Theorem 3 to develop an easily checked, necessary and sufficient condition for $\nu$-irreducibility of Eaton's Markov chain.

## 3. Conditions for irreducibility of Eaton's Markov chain

Eaton [3] studied the Markov transition function

$$R(\theta, \mathrm{d}\eta) = \int_{\mathsf{X}} Q(\mathrm{d}\eta|x) P(\mathrm{d}x|\theta),$$

where $P$ is the statistical model and $Q$ is an FPD. Since any FPD can be used to construct $R$, $R$ is not unique. For example, if $\tilde{Q} \neq Q$ is another FPD, then

$$\tilde{R}(\theta, \mathrm{d}\eta) = \int_{\mathsf{X}} \tilde{Q}(\mathrm{d}\eta|x) P(\mathrm{d}x|\theta)$$

is an equally valid version of the Markov transition function. However, the following result shows that $R$ enjoys a uniqueness property similar to the uniqueness property of the FPD that was discussed in Section 1.

**Proposition 2.** *If $R$ and $\tilde{R}$ are two different versions of Eaton's Markov transition function, then there exists a $\nu$-null set $B_0 \in \mathcal{B}(\Theta)$ such that $\theta \notin B_0$ implies $R(\theta, \cdot) = \tilde{R}(\theta, \cdot)$.*

**Proof.** Since $Q$ and $\tilde{Q}$ are both FPD's, there exists an $M$-null set $A_0$ such that $x \notin A_0 \Rightarrow Q(\mathrm{d}\theta|x) = \tilde{Q}(\mathrm{d}\theta|x)$. Fix $B \in \mathcal{B}(\Theta)$ and note that

$$R(\theta, B) - \tilde{R}(\theta, B) = \int_{A_0} [Q(B|x) - \tilde{Q}(B|x)] P(\mathrm{d}x|\theta).$$

Now, since $0 = M(A_0) = \int_\Theta P(A_0|\theta)\nu(\mathrm{d}\theta)$, there exists a $\nu$-null set $B_0$ such that $P(A_0|\theta) = 0$ for all $\theta \in \overline{B}_0$. Hence, $\theta \notin B_0 \Rightarrow R(\theta, B) = \tilde{R}(\theta, B)$. Finally, note that $A_0$ is determined by $Q$ and $\tilde{Q}$, and $B_0$ is determined by $A_0$. Therefore, $B_0$ does not depend on the set $B$, so $\theta \notin B_0 \Rightarrow R(\theta, \cdot) = \tilde{R}(\theta, \cdot)$. □

*Remark 1.* An important consequence of Proposition 2 is that either all versions of $R$ are locally-$\nu$-recurrent, or none of them is.

Every FPD satisfies $P(\mathrm{d}x|\theta)\nu(\mathrm{d}\theta) = Q(\mathrm{d}\theta|x)M(\mathrm{d}x)$ and it follows that every version of $R$ is reversible (or symmetric) with respect to $\nu$, that is, if $f$ and $g$ are bounded, real-valued functions on $\Theta$, then

$$\int_\Theta \int_\Theta f(\theta)g(\eta)R(\theta,\mathrm{d}\eta)\nu(\mathrm{d}\theta) = \int_\Theta \int_\Theta f(\theta)g(\eta)R(\eta,\mathrm{d}\theta)\nu(\mathrm{d}\eta).$$

This property is key in the proof of our main result, which we now state and prove.

**Theorem 4.** *There exists a $\nu$-reducible version of $R$ if and only if there exist a nonempty set $C \in \mathcal{B}(\Theta)$ with $\nu(\overline{C}) > 0$ and another set $A \in \mathcal{B}(\mathsf{X})$ such that $P(\overline{A}|\theta) = 0$ for every $\theta \in C$ and $P(A|\theta) = 0$ for $\nu$-a.e. $\theta \in \overline{C}$.*

**Proof.** To prove sufficiency, suppose that $A$ and $C$ exist. Using property 3 of Definition 1, we have

$$\int_A Q(\overline{C}|x)M(\mathrm{d}x) = \int_{\overline{C}} P(A|\theta)\nu(\mathrm{d}\theta) = 0.$$

Hence, if

$$D = \{x \in A : Q(\overline{C}|x) > 0\},$$

then $M(D) = 0$ and, obviously, $Q(\overline{C}|x) = 0$ for all $x \in A \setminus D$. Fix $\theta_0 \in C$ and let $\delta_{\theta_0}(\mathrm{d}\theta)$ denote a probability measure concentrated on the point $\theta_0$. Now, define

$$\tilde{Q}(\mathrm{d}\theta|x) = \begin{cases} Q(\mathrm{d}\theta|x), & \text{if } x \notin D, \\ \delta_{\theta_0}(\mathrm{d}\theta), & \text{if } x \in D. \end{cases}$$

Clearly, $\tilde{Q}$ satisfies the first and third properties of Definition 1. Moreover, we show in the Appendix that, for any $B \in \mathcal{B}(\Theta)$, $\tilde{Q}(B|\cdot)$ is measurable. Therefore, $\tilde{Q}$ is an FPD and we now show that

$$\tilde{R}(\theta, \mathrm{d}\eta) = \int_\mathsf{X} \tilde{Q}(\mathrm{d}\eta|x)P(\mathrm{d}x|\theta)$$

is the $\nu$-reducible version of $R$ that we seek. By construction, $\tilde{Q}(\overline{C}|x) = 0$ for all $x \in A$. It follows that, for every $\theta \in C$,

$$\tilde{R}(\theta, \overline{C}) = \int_\mathsf{X} \tilde{Q}(\overline{C}|x)P(\mathrm{d}x|\theta) = \int_A \tilde{Q}(\overline{C}|x)P(\mathrm{d}x|\theta) = 0.$$

Consequently, $C$ is a closed set with $\nu(\overline{C}) > 0$ and it follows from Theorem 3 that $\tilde{R}$ is $\nu$-reducible.

To prove necessity, assume that $R$ is $\nu$-reducible. By Theorem 3, there exists a closed set $C$ with $\nu(\overline{C}) > 0$. Using the reversibility of $R$, we have

$$\int_{\overline{C}} R(\theta, C)\nu(\mathrm{d}\theta) = \int_C R(\theta, \overline{C})\nu(\mathrm{d}\theta) = 0.$$

This, of course, implies that $R(\theta, C) = 0$ for $\nu$-a.e. $\theta \in \overline{C}$. Now, define $F_1 = \{x \in \mathsf{X} : Q(C|x) > 0\}$ and $F_2 = \{x \in \mathsf{X} : Q(\overline{C}|x) > 0\}$. Since $Q(C|x) + Q(\overline{C}|x) = 1$ for every $x \in \mathsf{X}$, we have $F_1 \cup F_2 = \mathsf{X}$. Consequently, $\overline{F}_2 \subset F_1$. Since $C$ is closed, we know that for any $\theta \in C$,

$$0 = R(\theta, \overline{C}) = \int_{\mathsf{X}} Q(\overline{C}|x) P(\mathrm{d}x|\theta) = \int_{F_2} Q(\overline{C}|x) P(\mathrm{d}x|\theta).$$

Thus, $P(F_2|\theta) = 0$ for all $\theta \in C$. Similarly, for $\nu$-a.e. $\theta \in \overline{C}$, we have

$$0 = R(\theta, C) = \int_{\mathsf{X}} Q(C|x) P(\mathrm{d}x|\theta) = \int_{F_1} Q(C|x) P(\mathrm{d}x|\theta).$$

Therefore, $P(F_1|\theta) = 0$ for $\nu$-a.e. $\theta \in \overline{C}$ and since $\overline{F}_2 \subset F_1$, it follows that $P(\overline{F}_2|\theta) = 0$ for $\nu$-a.e. $\theta \in \overline{C}$. Letting $A = \overline{F}_2$, we have $P(\overline{A}|\theta) = 0$ for all $\theta \in C$ and $P(A|\theta) = 0$ for $\nu$-a.e. $\theta \in \overline{C}$. □

Obviously, a sufficient condition for $\nu$-irreducibility of $R$ is the non-existence of the sets $A$ and $C$ in Theorem 4. This is precisely what Theorem 2 says. However, as will become clear from the examples in the next section, even when the sets $A$ and $C$ do exist, it is often the case that some versions of $R$ are $\nu$-irreducible. One way to establish local-$\nu$-recurrence (of all versions) of $R$ is to identify a single version of $R$ that is $\nu$-irreducible and then show that this version of $R$ is recurrent. The result then follows from the equivalence of recurrence and local-$\nu$-recurrence (under irreducibility) mentioned in Section 1.

## 4. Examples

In this section, we illustrate the use of Theorem 4 with two examples.

**Example 1.** Let $\mathsf{X} = \Theta = \mathbb{R}$ and $P(\mathrm{d}x|\theta) = p(x|\theta)\,\mathrm{d}x$, where

$$p(x|\theta) = I_{(\theta, \theta+1)}(x)$$

and where $\mathrm{d}x$ denotes Lebesgue measure on $\mathsf{X}$. Take the prior distribution to be $\nu(\mathrm{d}\theta) = \mathrm{d}\theta$. A simple calculation shows that the marginal measure is $M(\mathrm{d}x) = \mathrm{d}x$, which is clearly $\sigma$-finite. This is all the information we require to apply Theorem 4. Let $C \in \mathcal{B}(\Theta)$ be a non-empty set such that $\nu(\overline{C}) > 0$. We claim that there exist $\theta' \in C$ and $D \subset \overline{C}$ such that $\nu(D) > 0$ and $|\theta - \theta'| < 1/4$ for every $\theta \in D$. To see this, let $D_i = [\frac{i}{8}, \frac{i+1}{8})$ for all $i \in \mathbb{Z} := \{\ldots, -1, 0, 1, \ldots\}$. Then, $\Theta = \bigcup_{i \in \mathbb{Z}} D_i$. If there exists an $i_0$ such that $C \cap D_{i_0}$ is not empty and $\nu(\overline{C} \cap D_{i_0}) > 0$, then we can simply take $\theta'$ to be any point in $C \cap D_{i_0}$ and $D$ to be the set $\overline{C} \cap D_{i_0}$. Otherwise, for every $i \in \mathbb{Z}$, either $D_i \subset \overline{C}$ or $\nu(\overline{C} \cap D_i) = 0$. Since $\nu(\overline{C}) > 0$, there must exist an $i_0$ such that $\nu(\overline{C} \cap D_{i_0}) > 0$ and it follows that $D_{i_0} \subset \overline{C}$. Assume, without loss of generality, that $C \cap [\frac{i_0}{8}, \infty)$ is non-empty and let $i_1 = \min\{i > i_0 : C \cap D_i \text{is non-empty}\}$, which is clearly finite. Now, any point in $C \cap D_{i_1}$ and the set $D_{i_1-1}$ play the roles of $\theta'$ and $D$. These arguments show that there exist $\theta \in \overline{C}$ and $\theta' \in C$ such that $|\theta - \theta'| < 1/4$, where $\theta$ can be chosen outside of any subset of $\overline{C}$ having Lebesgue measure zero. It follows that if $I_1 = (\theta, \theta+1)$, $I_2 = (\theta', \theta'+1)$ and $I = I_1 \cap I_2$, then $\nu(I) > 0$. Now, if there exists an $A \in \mathcal{B}(\mathsf{X})$ such that $P(A|\theta) = 0$ and $P(\overline{A}|\theta') = 0$, then

$$\nu(I) = \nu(A \cap I) + \nu(\overline{A} \cap I) \le \nu(A \cap I_1) + \nu(\overline{A} \cap I_2) = 0,$$

which is a contradiction. Hence, such an $A$ cannot exist and it follows from Theorem 4 that every version of Eaton's Markov chain is $\nu$-irreducible.

For the sake of comparison, we now explain what is required to establish irreducibility in this situation if we make no appeal to Theorem 4. A version of the posterior is given by $Q(\mathrm{d}\theta|x) = q(\theta|x)\,\mathrm{d}\theta$, where

$$q(\theta|x) = I_{(x-1,x)}(\theta).$$

It follows that $R(\theta, \mathrm{d}\eta) = r(\eta|\theta)\,\mathrm{d}\eta$, where

$$r(\eta|\theta) = (1 + (\eta - \theta))I_{(-1,0)}(\eta - \theta) + (1 - (\eta - \theta))I_{(0,1)}(\eta - \theta).$$

Since $\theta$ is a location parameter in the density $r(\eta|\theta)$, the Markov chain $W$ can be expressed as a random walk

$$W_{n+1} = W_n + Z_{n+1},$$

where $Z_1, Z_2, \ldots$ is an independent and identically distributed (i.i.d.) sequence of random variables with density given by

$$f(z) = (1+z)I_{(-1,0)}(z) + (1-z)I_{(0,1)}(z).$$

For any $\varepsilon \in (0,1)$, we have

$$P(Z_1 \in (0, \varepsilon)) = P(Z_1 \in (-\varepsilon, 0)) > 0.$$

This implies that the chain $W$ can make arbitrarily small jumps in either direction. While this argument makes it intuitively clear that the chain is $\nu$-irreducible, a formal proof requires a technical argument similar to that used in Section 4.3.3 of Meyn and Tweedie [11].

It turns out that this random walk is recurrent, which implies that $W$ is locally-$\nu$-recurrent. Hence, by Theorem 1, $\nu$ is strongly admissible. In fact, Eaton [3] shows that, under very mild conditions, Lebesgue measure is a strongly admissible prior for one- and two-dimensional location problems. This concludes Example 1.

***Example 2.*** Let $\Theta = \mathbb{R}_+ := [0, \infty)$ and $\mathsf{X} = \mathbb{R}_+^n$. Suppose that when $\theta > 0$, the statistical model $P(\mathrm{d}x|\theta)$ has a density (with respect to Lebesgue measure on $\mathsf{X}$) given by

$$p(x|\theta) = \prod_{i=1}^n \frac{1}{\theta} I_{[0,\theta)}(x_i),$$

where $x = (x_1, \ldots, x_n)$. In words, our statistical model stipulates that $X_1, \ldots, X_n$ are i.i.d. random variables from the uniform distribution on $[0, \theta)$. Take the prior distribution to be $\nu(\mathrm{d}\theta) = \frac{\mathrm{d}\theta}{\theta}$, where $\mathrm{d}\theta$ denotes Lebesgue measure on $\Theta$. While we have yet to define $P(\mathrm{d}x|\theta)$ in the case where $\theta = 0$, from a practical (statistical) standpoint this definition is irrelevant since $\nu(\{0\}) = 0$. However, technically speaking, this distribution must be specified to complete the model. (Note that we cannot simply remove the point $\{0\}$ from $\Theta$ because $(0, \infty)$ is not a Polish space.) We consider two alternatives for $P(\mathrm{d}x|0)$:

1. a unit point mass at $(0, \ldots, 0) \in \mathbb{R}^n$;
2. $n$ i.i.d. exponential random variables with unit scale.

We now apply Theorem 4. In case 1, there does exist a $\nu$-reducible version of $R$ since we can take $A$ to be the point $(0, \ldots, 0) \in \mathbb{R}^n$ and $C = \{0\}$. However, in case 2, the sets $A$ and $C$ do not exist. Indeed, let $\mu$ denote Lebesgue measure on $\mathsf{X}$ and note that it is impossible to satisfy the conditions of Theorem 4 if $\mu(\overline{A}) = 0$ since this yields $P(A|\theta) = 1$ for all $\theta \in \Theta$. Now, if $\mu(\overline{A}) > 0$ and $P(\overline{A}|\theta) = 0$, then it must be true that $\theta > 0$ and $\mu(\overline{A} \cap [0, \theta)^n) = 0$. Similarly, if $\theta^* > 0$ and $P(A|\theta^*) = 0$, then $\mu(A \cap [0, \theta^*)^n) = 0$. Now, let $\theta_0 = \min\{\theta, \theta^*\}$ and note that

$$\mu([0, \theta_0)^n) = \mu(A \cap [0, \theta_0)^n) + \mu(\overline{A} \cap [0, \theta_0)^n) \leq \mu(A \cap [0, \theta^*)^n) + \mu(\overline{A} \cap [0, \theta)^n) = 0.$$

This is a contradiction, which implies that the sets $A$ and $C$ do not exist. It follows that every version of $R$ is $\nu$-irreducible. Note that, by defining the statistical model carefully on an irrelevant ($\nu$-null) set of $\theta$'s, we were able to employ Theorem 4 to show that all versions of $R$ are $\nu$-irreducible.

Regardless of how $P(\mathrm{d}x|0)$ is defined, the marginal measure is given by

$$M(\mathrm{d}x) = \frac{\mathrm{d}x}{nx_{(n)}^n},$$

where $x_{(n)} := \max\{x_1, \ldots, x_n\}$ and $\mathrm{d}x$ denotes Lebesgue measure on $\mathsf{X}$. Consider case 1 again. For $x$ such that $x_{(n)} > 0$, define

$$q(\theta|x) = \frac{nx_{(n)}^n}{\theta^{n+1}} I_{(x_{(n)}, \infty)}(\theta).$$

Since the point $(0, \ldots, 0) \in \mathbb{R}^n$ has $M$-measure 0, $Q(\mathrm{d}\theta|x)$ can essentially be chosen arbitrarily when $x_{(n)} = 0$. We consider two different choices. Let $\delta_0(\mathrm{d}\theta)$ denote the probability measure concentrated at $\{0\}$ and let $\delta_1(\mathrm{d}\theta)$ denote a probability measure with support $\mathbb{R}_+$. (It seems more appropriate to take $Q(\mathrm{d}\theta|x)$ equal to $\delta_0(\mathrm{d}\theta)$ when $x_{(n)} = 0$.) Two different versions of the posterior distribution are given by

$$Q_i(\mathrm{d}\theta|x) = \begin{cases} q(\theta|x)\,\mathrm{d}\theta, & x_{(n)} > 0, \\ \delta_i(\mathrm{d}\theta), & x_{(n)} = 0, \end{cases}$$

for $i \in \{0,1\}$. The version of Eaton's chain associated with $Q_0$ is $\nu$-reducible since, if the chain is started at $\theta = 0$, it stays there forever. On the other hand, Eaton *et al.* [7] show that the version associated with $Q_1$ is $\nu$-irreducible and go on to show that this chain is recurrent. Consequently, this version of Eaton's chain is locally-$\nu$-recurrent and it follows from Theorem 1 that $\nu$ is strongly admissible. This concludes Example 2.

## Appendix: On the measurability of $\tilde{Q}$

Here, we establish the measurability of $\tilde{Q}$. Fix $B \in \mathcal{B}(\Theta)$. It suffices to show that, for any $t \in \mathbb{R}$, the set

$$\tilde{G}_t = \{x \in \mathsf{X} : \tilde{Q}(B|x) < t\}$$

is in $\mathcal{B}(\mathsf{X})$. We will accomplish this using the partition

$$\tilde{G}_t = (\tilde{G}_t \cap D) \cup (\tilde{G}_t \cap \overline{D}).$$

The set $G_t = \{x \in \mathsf{X} : Q(B|x) < t\} \in \mathcal{B}(\mathsf{X})$ since $Q(B|\cdot)$ is a measurable function. The measures $Q$ and $\tilde{Q}$ agree on $\overline{D}$, hence

$$\tilde{G}_t \cap \overline{D} = \{x \in \overline{D} : \tilde{Q}(B|x) < t\} = \{x \in \overline{D} : Q(B|x) < t\} = G_t \cap \overline{D} \in \mathcal{B}(\mathsf{X}).$$

It remains to show that $\tilde{G}_t \cap D \in \mathcal{B}(\mathsf{X})$. There are four possible cases:

1. $\theta_0 \in B$ and $t \leq 1$;
2. $\theta_0 \in B$ and $t > 1$;
3. $\theta_0 \notin B$ and $t \leq 0$;
4. $\theta_0 \notin B$ and $t > 0$.

In the first two cases, $\tilde{Q}(B|x) = 1$ for all $x \in D$. Hence, in the first case, $\tilde{G}_t \cap D = \emptyset \in \mathcal{B}(\mathsf{X})$, and in the second, $\tilde{G}_t \cap D = D \in \mathcal{B}(\mathsf{X})$. In the last two cases, $\tilde{Q}(B|x) = 0$ for all $x \in D$. Hence, in the third case, $\tilde{G}_t \cap D = \emptyset \in \mathcal{B}(\mathsf{X})$, and in the last case, $\tilde{G}_t \cap D = D \in \mathcal{B}(\mathsf{X})$.

## Acknowledgements

The authors thank three anonymous reviewers for helpful comments and suggestions. Hobert's research was partially supported by NSF Grant DMS-05-03648.